\documentclass[12pt,leqno,draft]{article}

\newtheorem{theorem}{Theorem}
\newtheorem{lemma}[theorem]{Lemma}
\newtheorem{proposition}[theorem]{Proposition}
\newtheorem{definition}[theorem]{Definition}
\newtheorem{corollary}[theorem]{Corollary}

\newcommand{\begintheorem}{\addtocounter{equation}{1}\begin{theorem}}
\newcommand{\beginlemma}{\addtocounter{equation}{1}\begin{lemma}}
\newcommand{\beginproposition}{\addtocounter{equation}{1}\begin{proposition}}
\newcommand{\begindefinition}{\addtocounter{equation}{1}\begin{definition}}
\newcommand{\begincorollary}{\addtocounter{equation}{1}\begin{corollary}}

\begin{document}

\title{Some notes about matrices, 3}

\author{Stephen William Semmes	\\
	Rice University		\\
	Houston, Texas}

\date{}

\maketitle

\tableofcontents

\section{Projective spaces}
\label{section on projective spaces}
\setcounter{equation}{0}

	Let $n$ be a positive integer.  As usual, ${\bf R}$, ${\bf C}$
denote the real and complex numbers, respectively, and ${\bf R}^n$,
${\bf C}^n$ consist of the $n$-tuples of real and complex numbers,
respectively.  The $n$-dimensional \emph{real and complex projective
spaces} ${\bf RP}^n$, ${\bf CP}^n$ consist of the real and complex lines
through the origin in ${\bf R}^{n+1}$, ${\bf C}^{n+1}$, respectively.

	To put it another way, if ${\bf R}_*$, ${\bf C}_*$ denote the
nonzero real and complex numbers, respectively, then we have natural
actions of ${\bf R}_*$, ${\bf C}_*$ on ${\bf R}^{n+1} \backslash
\{0\}$, ${\bf C}^{n+1} \backslash \{0\}$ by scalar multiplication, and
the projective spaces are the corresponding quotient spaces.  In other
words, two nonzero vectors $v$, $w$ in ${\bf R}^{n+1}$, ${\bf
C}^{n+1}$ lead to the same point in the corresponding projective space
exactly when they are scalar multiples of each other.  Note that we
get canonical mappings from ${\bf R}^{n+1} \backslash \{0\}$, ${\bf
C}^{n+1} \backslash \{0\}$ onto ${\bf RP}^n$, ${\bf CP}^n$, in which a
nonzero vector $v$ is sent to the line through the origin which passes
through $v$, which consists of all scalar multiples of $v$.

	If $L$ is a nontrivial linear subspace of ${\bf R}^{n+1}$,
${\bf C}^{n+1}$ of dimension $l + 1$, say, then we get an interesting
space ${\bf P}(L)$ consisting of all lines through the origin in $L$,
which we can think of as sitting inside of ${\bf RP}^n$, ${\bf CP}^n$,
as appropriate.  More precisely, ${\bf P}(L)$ is basically a copy of
${\bf RP}^l$ or ${\bf CP}^l$.  These are the $l$-dimensional ``linear
subspaces'' of projective space, analogous to linear subspaces of
${\bf R}^n$, ${\bf C}^n$.

	If $A$ is an invertible linear transformation on ${\bf
R}^{n+1}$ or on ${\bf C}^{n+1}$, then $A$ takes lines to lines, and
induces a transformation $\widehat{A}$ on the corresponding projective
space.  Notice that $\widehat{A}$ is automatically a one-to-one
transformation of the corresponding projective space onto itself, with
\begin{equation}
	\widehat{A}^{-1} = \widehat{(A^{-1})},
\end{equation}
and $\widehat{A}$ maps linear subspaces of projective space to
themselves, in the sense of the preceding paragraph.  Also, if $A_1$,
$A_2$ are invertible linear transformations on ${\bf R}^{n+1}$ or on
${\bf C}^{n+1}$, then the induced transformations $\widehat{A}_1$,
$\widehat{A}_2$ on the corresponding projective space satisfy
\begin{equation}
	\widehat{A}_1 \circ \widehat{A}_2 = \widehat{(A_1 \circ A_2)}.
\end{equation}

	Let $H$ be a hyperplane in ${\bf R}^{n+1}$ or in ${\bf C}^{n+1}$,
which is to say a linear subspace of dimension $n$, and let $v$ be a
nonzero vector in ${\bf R}^{n+1}$, ${\bf C}^{n+1}$, as appropriate.
This leads to an affine hyperplane $H + v$, consisting of all vectors
of the form $w + v$, $w \in H$, and which does not contain the vector
$0$.  For each $w \in H$, we can look at the line through $w + v$,
which we can view as an element of the corresponding projective space.

	In other words, we basically get an embedding of $H$ into the
corresponding projective space, ${\bf RP}^n$ or ${\bf CP}^n$.  Of
course we can also think of $H$ as being isomorphic to ${\bf R}^n$ or
${\bf C}^n$, so that we are really looking at a bunch of embeddings of
${\bf R}^n$, ${\bf C}^n$ into ${\bf RP}^n$, ${\bf CP}^n$,
respectively.  For instance, we can do this with $H$ equal to the
$j$th coordinate hyperplane in ${\bf R}^{n+1}$, ${\bf C}^{n+1}$,
$1 \le j \le n + 1$, which is defined by the condition that the
$j$th coordinate of vectors in $H$ are equal to $0$, and we can take
$v$ to be the $j$th standard basis vector, with $j$th coordinate
equal to $1$ and the other $n$ coordinates equal to $0$.

	These $n + 1$ embeddings of ${\bf R}^n$, ${\bf C}^n$ into
${\bf RP}^n$, ${\bf CP}^n$ corresponding to the $n + 1$ coordinate
hyperplanes in ${\bf R}^{n+1}$, ${\bf CP}^{n+1}$ are sufficient to
cover the projective space, i.e., every point in projective space
shows up in the image of at least one of the embeddings.  For a given
hyperplane $H$, the set of points in the projective space which do not
occur in the embedding of $H$ is the same as ${\bf P}(H)$.  Thus
the set of missing points in the projective space lie in a projective
subspace of dimension $1$ less.

	Using these embeddings of ${\bf R}^n$, ${\bf C}^n$ into
the corresponding projective spaces, we can think of the projective spaces
as being manifolds.  That is, these embeddings provide local coordinates
for all points in the projective space.  Two different embeddings
which contain the same point $p$ in their image are compatible in terms
of topology and also smooth structure, and indeed there is a finer
``projective'' structure which is reflected in the presence of nice
projective subspaces, for instance.

	Note that two invertible linear transformations $A_1$, $A_2$
on ${\bf R}^{n+1}$ or on ${\bf C}^{n+1}$ lead to the same induced
transformation on projective space if and only if there is a nonzero
scalar $\alpha$ such that $A_2 = \alpha \, A_1$.  Thus the group of
these ``projective linear transformations'' has dimension $(n+1)^2 -
1$ over the real or complex numbers, as appropriate.  Also, for
any pair of points $p$, $q$ in a projective space, there is a
projective linear transformation which takes $p$ to $q$.

\section{Grassmannians}
\label{section on grassmannians}
\setcounter{equation}{0}

	Fix positive integers $k$, $n$ with $k < n$.  The Grassmann
spaces $G_{\bf R}(k, n)$, $G_{\bf C}(k, n)$ consist of the
$k$-dimensional linear subspaces of ${\bf R}^n$, ${\bf C}^n$,
respectively.  When $k = 1$ this reduces to the $(n-1)$-dimensional
projective spaces.

	Suppose that $L$, $M$ are linear subspaces of ${\bf R}^n$
or of ${\bf C}^n$ such that $L$ has dimension $k$, $M$ has dimension
$n - k$, and the intersection of $L$, $M$ consists of only the zero
vector.  Thus $L$ and $M$ are complementary, and if $A$ is a linear
mapping from $L$ to $M$, then the graph of $A$, consisting of the
vectors
\begin{equation}
	v + A(v), \quad v \in L,
\end{equation}
is also a $k$ dimensional subspace of ${\bf R}^n$ or of ${\bf C}^n$,
as appropriate.  In this way we can embed the vector space of linear
transformations from $L$ to $M$ into the Grassmannian, and in
particular this provides a nice coordinate patch around $L$ itself.

	In particular, notice that the dimension of the Grassmann
space of $k$ planes in ${\bf R}^n$ or ${\bf C}^n$ is
\begin{equation}
	k (n-k)
\end{equation}
with respect to the real or complex numbers, as appropriate.  Just as
for projective spaces, invertible linear transformations on ${\bf
R}^n$ or on ${\bf C}^n$ induce interesting mappings on the
corresponding Grassmannians.  These actions are again transitive,
because if $L_1$, $L_2$ are $k$-dimensional linear subspaces of ${\bf
R}^n$ or of ${\bf C}^n$, then there is an invertible linear
transformation $A$ on ${\bf R}^n$ or on ${\bf C}^n$, as appropriate,
such that $A(L_1) = L_2$.

	There is a natural correspondence between the Grassmann spaces
of $k$-dimensional linear subspaces in ${\bf R}^n$ or in ${\bf C}^n$
and the Grassmann spaces of $(n-k)$-dimensional linear subspaces of
${\bf R}^n$ or in ${\bf C}^n$, respectively.  One might prefer to
think of $k$-dimensional linear subspaces of ${\bf R}^n$ or ${\bf
C}^n$ as being associated to $(n-k)$-dimensional linear subspaces of
the corresponding dual spaces, by looking at intersections of kernels
of linear functionals in a subspace of a dual space.  Alternatively,
one can think of linear subspaces of dimension $k$ as being associated
to their orthogonal complements, which are linear subspaces of
dimension $n - k$, using inner products.

\section{Hopf-type spaces}
\label{section on Hopf-type spaces}
\setcounter{equation}{0}

	Fix a positive integer $n$, and consider ${\bf R}^n \backslash
\{0\}$, ${\bf C}^n \backslash \{0\}$.  On these spaces one can consider
multiplication by integer powers of $2$.  One can then consider the
corresponding quotients by this action, which is to say to identify
two nonzero vectors which can be expressed as integer powers of $2$
times each other.

	There are natural mappings from these spaces to the
corresponding projective spaces of dimension $n - 1$, by looking at
the lines through nonzero vectors, which are not changed by nonzero
scalar multiples in general.  As for projective spaces, invertible
linear mappings on ${\bf R}^n$ or on ${\bf C}^n$ induce nice
transformations on the quotients by the actions of powers of
$2$, and indeed these transformations are compatible with the
ones on projective spaces too.  In the quotients by integer
powers of $2$, there are interesting subspaces arising from linear
subspaces back in ${\bf R}^n$ or in ${\bf C}^n$.

	Of course one can consider integer powers of other scalars
with absolute value larger than $1$ in place of $2$, and have similar
properties.  Analogous objects have been studied using more
complicated invertible linear transformations on ${\bf R}^n$ or on
${\bf C}^n$ rather than scalar multiplication.  This can lead to more
tricky kinds of twisting, for instance.

\section{Hopf fibrations}
\label{section on Hopf fibrations}
\setcounter{equation}{0}

	Let $n$ be a positive integer.  There is an obvious mapping
from the unit sphere in ${\bf R}^{n+1}$, consisting of vectors with
norm $1$, to $n$-dimensional real projective space ${\bf RP}^n$.
Namely, if $x \in {\bf R}^{n+1}$ and $|x| = 1$, the line through $x$
leads to a point in ${\bf RP}^n$, and every element of ${\bf RP}^n$
arises in this manner exactly twice.

	Similarly, there is a natural mapping from the unit sphere in
${\bf C}^{n+1}$, consisting of vectors with norm $1$, to complex
projective space ${\bf CP}^n$.  If $v$ is a vector in ${\bf C}^{n+1}$
with $|v| = 1$, then the line through $v$ leads to a point in ${\bf
C}^n$.  Two vectors $v$, $w$ in ${\bf C}^n$ with $|v| = |w| = 1$
correspond to the same element of ${\bf CP}^n$ if and only if there is
a complex number $\alpha$ such that $|\alpha| = 1$ and $w = \alpha \,
v$.

	In other words, the inverse image of a point in ${\bf RP}^n$
back in the unit sphere of ${\bf R}^{n+1}$ consists of two points,
while in the complex case the inverse image of a point in ${\bf CP}^n$
back in the unit sphere in ${\bf C}^{n+1}$ is a circle.  In terms
of the unit sphere in ${\bf C}^{n+1}$, there is a family of circles
which cover the whole sphere and which are pairwise disjoint.
When $n = 1$, these are circles in the $3$-dimensional sphere
which are linked with each other.

	The $n = 1$ case has other nice properties.  One can identify
${\bf CP}^1$ with a $2$-dimensional sphere, basically because ${\bf
CP}^1$ can be expressed as ${\bf C}$ with a single additional point,
which one might denote $\infty$.  With respect to this identification,
projective linear transformations on ${\bf CP}^1$ correspond to
transformations of the form
\begin{equation}
	z \mapsto \frac{a \, z + b}{c \, z + d}
\end{equation}
on ${\bf C} \cup \{\infty\}$, with the usual conventions like $1/0 =
\infty$, and where $({a b \atop c d})$ is an invertible $2 \times 2$
matrix.

\end{document}